\documentclass[12pt]{article}

\usepackage{amsmath,amssymb,amscd}
\usepackage{graphicx, psfrag}

\begin{document}

\title{Structure of Three-Manifolds\\
-- Poincar\'{e} and geometrization conjectures}
\author{Shing-Tung Yau$^{1,2}$}
\date{}

\maketitle

\footnotetext[1]{This was a talk given at the Morningside Center
of Mathematics on June 20, 2006.}

\footnotetext[2]{All the computer graphics are provided by David
Gu, based on the joint paper of David Gu, Yalin Wang and S.-T.
Yau.}

Ladies and gentlemen, today I am going to tell you the story of
how a chapter of mathematics has been closed and a new chapter is
beginning.

Let me begin with some elementary observations.

A major purpose of Geometry is to describe and classify geometric
structures of interest. We see many such interesting structures in
our day-to-day life.

Let us begin with topological structures of a two dimensional
surface. These are spaces where locally we have two degrees of
freedom. Here are some examples:

\begin{figure}[ht]
\begin{center}
\begin{tabular}{cccc}
\includegraphics[height=2.25cm]{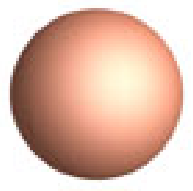}&
\includegraphics[height=2.25cm]{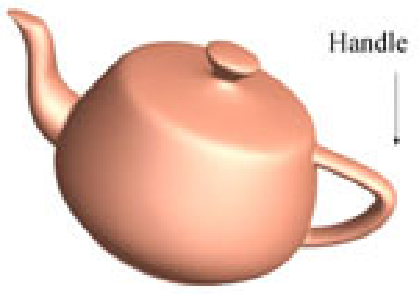}&
\includegraphics[height=2.25cm]{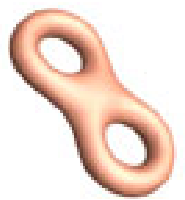}&
\includegraphics[height=2.25cm]{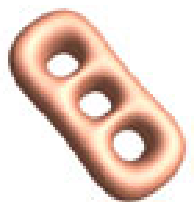}\\
genus 0& genus 1& genus 2& genus 3
\end{tabular}
\end{center}
\end{figure}

Genus of a surface is the number of handles of the surface.

An abstract and major way to construct surfaces is by connecting
along some deleted disk of each surface.

The connected sum of two surfaces $S_1$ and $S_2$ is denoted by
$S_1\# S_2$. It is formed by deleting the interior of disks $D_i$
from each $S_i$ and attaching the resulting punctured surfaces
$S_i-D_i$ to each other by a one-to-one continuous map $h:
\partial D_1 \to
\partial D_2$, so that
\[
    S_1 \# S_2 = (S_1 - D_1) \cup_h (S_2 - D_2).
\]

\begin{figure}[ht]
\psfrag{M1}{$S_1$} \psfrag{M2}{$S_2$}
\psfrag{B1}{$D_1$}\psfrag{B2}{$D_2$}
\begin{center}
\begin{tabular}{c}
\includegraphics[height=2.0cm]{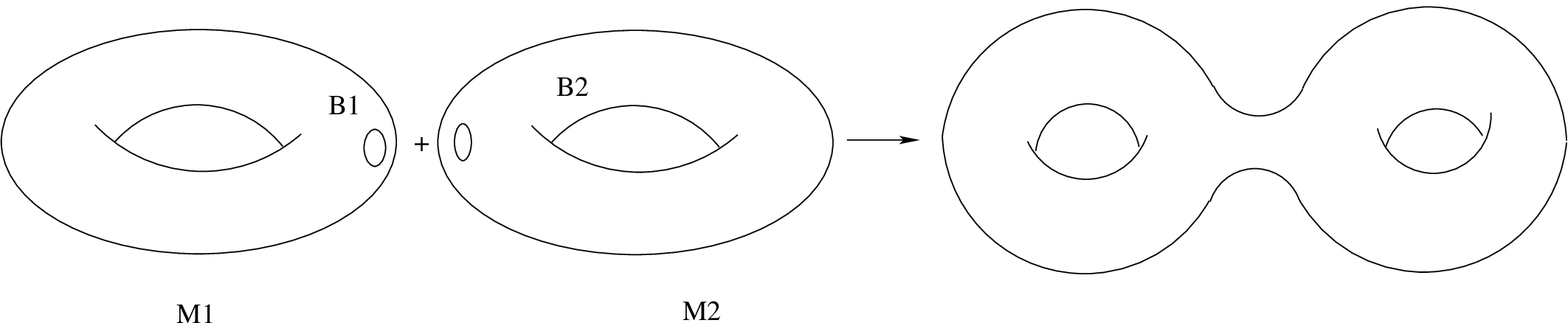}
\end{tabular}
\end{center}
\end{figure}

Example:

\begin{figure}[ht]
\begin{center}
\begin{tabular}{c}
\includegraphics[height=5.5cm]{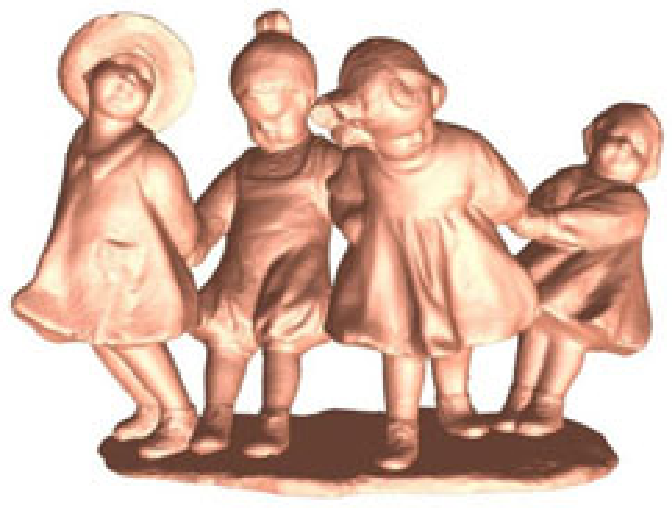}
\end{tabular}\\
A genus 8 surface, constructed by connected sum.
\end{center}
\end{figure}

The major theorem for the two dimensional surfaces is the
following:

\medskip

\noindent{\bf Theorem (Classification Theorem for Surfaces).} {\em
Any closed, connected orientable surface is exactly one of the
following surfaces: a sphere, a torus, or a finite number of
connected sum of tori. }

Note that a surface is called orientable if each closed curve on
it has a well-defined continuous normal field.

\section{Conformal geometry}

In order to understand surfaces in a deep manner, Riemann,
Poincar\'e  and others proposed to study conformal structure on
these two dimensional objects. Such structures allow us to measure
angles in the neighborhood of each point on the surface.

For example, if we take a standard atlas of the globe, we have
longitude and latitude. They are orthogonal to each other. When we
map the atlas, which is a square, onto the globe; distances are
badly distorted. For example, the region around the north pole is
shown to be a large region on the square. However, the fact that
longitude and latitude is orthogonal to each other is preserved
under the map. Hence if a ship moves in the ocean, we can use the
atlas to determine its direction accurately, but not the distance
travelled.

Globe

\begin{figure}[ht]
\begin{center}
\begin{tabular}{cc}
\includegraphics[height=3.0cm]{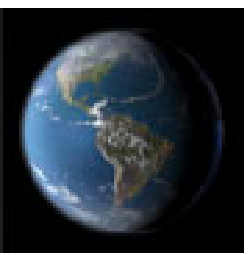}&
\includegraphics[height=2.75cm]{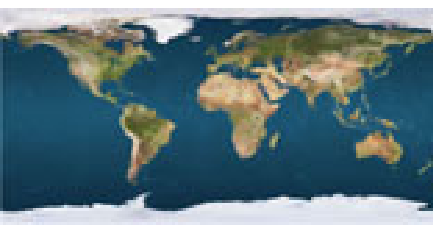}\\
\end{tabular}
\end{center}
\end{figure}

Poincar\'e found that at any point, we can draw a longitude (blue
curve) and latitude (red curve) on any surface of genus zero in
three space. These curves are orthogonal to each other and they
converge to two distinct points, on the surface, just like north
pole and south pole on the sphere. This theorem of Poincar\'{e}
also works for arbitrary abstract surface with genus zero.

It is a remarkable theorem that for any two closed surfaces with
genus zero, we can always find a one-to-one continuous map mapping
longitude and latitude of one surface to the corresponding
longitude and latitude of the other surface. This map preserve
angles defined by the charts. In such a situation, we say that
these two surfaces are conformal to each other. And there is only
one conformal structure for a surface with genus zero.

For genus equal to one, the surface looks like a donut, and we can
draw longitude and latitude with no north or south poles. However,
there can be distinct surfaces with genus one that are not
conformal to each other. In fact, there are two parameters of
conformal structures on a genus one surface. For genus $g$ greater
than one, one can still draw longitude and latitude (the
definition of such curves needs to be made precise). But they have
many poles, the number of which depends on the genus. The number
of parameters of conformal structures over a surface with genus
$g$ is $6g-6$.

In order to find a global atlas of the surface, we can cut along
some special curves of a surface and then spread the surface on
the plane or the disk. In this procedure, the longitude and the
latitude will be preserved.
\begin{figure}[ht]
\begin{center}
\begin{tabular}{ccc}
\includegraphics[height=2.75cm]{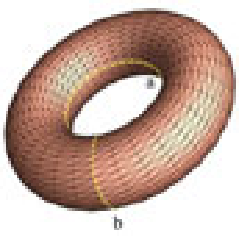}&
\includegraphics[height=2.75cm]{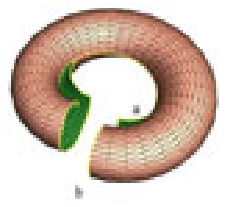}&
\includegraphics[height=2.75cm]{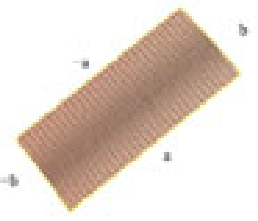}\\
\end{tabular}
\begin{tabular}{cc}
\includegraphics[height=3.0cm]{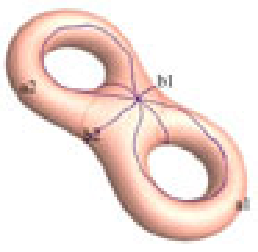}&
\includegraphics[height=3.0cm]{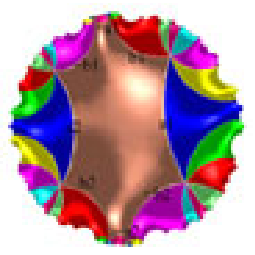}\\
\end{tabular}
\end{center}
\end{figure}

A fundamental theorem for surfaces with metric structure is the
following theorem.

\medskip

\noindent\textbf{Theorem (Poincar\'{e}'s Uniformization Theorem).}
{\em Any closed two-dimensional space
 is conformal to another space with constant Gauss curvature.
 \begin{itemize}
    \item If curvature $>0$, the surface has genus $=0$;
    \item If curvature $=0$, the surface has genus $=1$;
    \item If curvature $<0$, the surface has genus $>1$.
 \end{itemize} }

The generalization of this theorem plays a very important role in
the field of geometric analysis. In particular, it motivates the
works of Thurston and Hamilton. This will be discussed later in
this talk.
\begin{figure}[ht]
\begin{center}
\begin{tabular}{ccc}
\includegraphics[height=2.5cm]{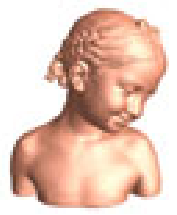}&
\includegraphics[height=2.5cm]{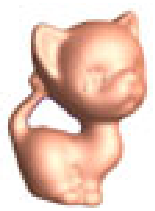}&
\includegraphics[height=2.5cm]{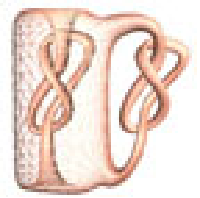}\\
\includegraphics[height=2.5cm]{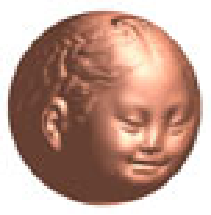}&
\includegraphics[height=2.5cm]{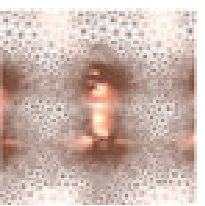}&
\includegraphics[height=2.5cm]{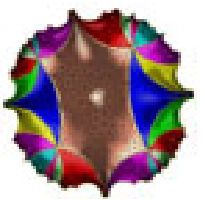}\\
Spherical & Euclidean & Hyperbolic \\
\end{tabular}
\end{center}
\end{figure}

\section{Hamilton's equation on Surfaces}

Poincar\'{e}'s theorem can also be proved by the equation of
Hamilton. We can deform any metric on a surface by the negative of
its curvature. After normalization, the final state of such
deformation will be a metric with constant curvature. This is a
method created by Hamilton to deform metrics on spaces of
arbitrary dimensions. In higher dimension, the typical final state
of spaces for the Hamilton equation is a space that satisfies
Einstein's equation.

As a consequence of the works by Richard Hamilton and B. Chow, one
knows that in two dimension, the deformation encounters no
obstruction and will always converge to one with constant
curvature. This theorem was used by David Gu, Yalin Wang, and
myself for computer graphics. The following sequence of pictures
is obtained by numerical simulation of the Ricci flow in two
dimension.
\begin{figure}[ht]
\begin{center}
\begin{tabular}{cccc}
\includegraphics[height=2cm]{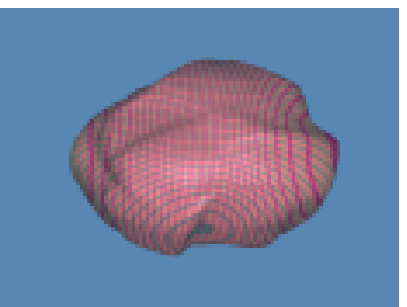}&
\includegraphics[height=2cm]{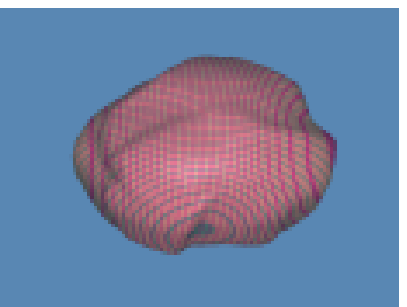}&
\includegraphics[height=2cm]{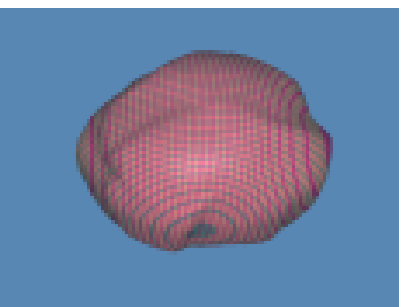}&
\includegraphics[height=2cm]{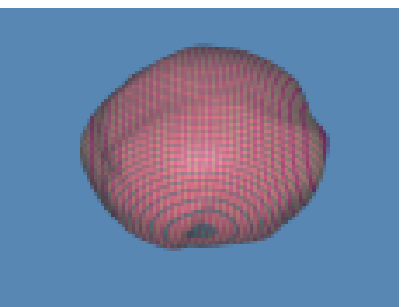}\\
\includegraphics[height=2cm]{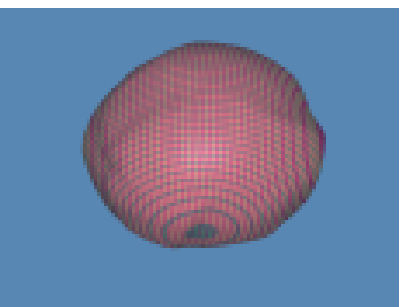}&
\includegraphics[height=2cm]{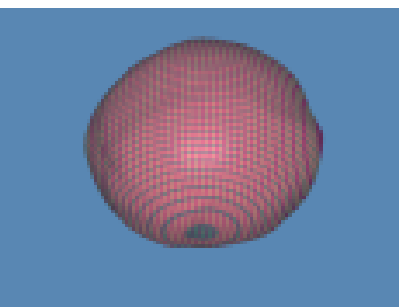}&
\includegraphics[height=2cm]{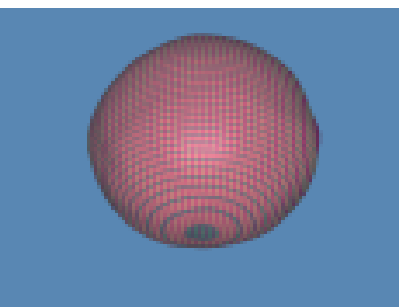}&
\includegraphics[height=2cm]{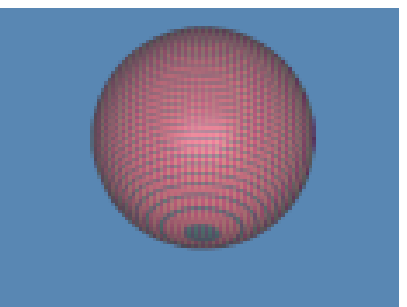}\\
\end{tabular}
\end{center}
\end{figure}

\section{Three-Manifolds}

So far, we have focused on spaces where there are only two degrees
of freedom. Instead of being a flat bug moving with two degrees of
freedom on a surface, we experience three degrees of freedom in
space. While it seems that our three dimensional space is flat,
there are many natural three dimensional spaces, which are not
flat.

Important natural example of higher dimensional spaces are phase
spaces in mechanics.

In the early twentieth century, Poincar\'{e} studied the topology
of phase space of dynamics of particles. The phase space consists
of $(x; v)$, the position and the velocity of the particles. For
example if a particle is moving freely with unit speed on a two
dimensional surface $\Sigma$, there are three degrees of freedom
in the phase space of the particle. This gives rise to a three
dimensional space $M$.

Such a phase space is a good example for the concept of fiber
bundle.

If we associate to each point $(x; v)$ in $M$ the point
$x\in\Sigma$, we have a map from $M$ onto $\Sigma$. When we fix
the point $x$, $v$ can be any vector with unit length. The
totality of $v$ forms a circle. Therefore, $M$ is a fiber bundle
over $\Sigma$ with fiber equal to a circle.

\section{The Poincar\'{e} Conjecture}

The subject of higher dimensional topology started with
Poincar\'{e}'s question:

Is a closed three dimensional space topologically a sphere if
every closed curve in this space can be shrunk continuously to a
point?

This is not only a famous difficult problem, but also the central
problem for three dimensional topology. Its understanding leads to
the full structure theorem for three dimensional spaces. I shall
describe its development chronologically.

\section{Topological Surgery}

Topologists have been working on this problem for over a century.
The major tool is application of cut and paste, or surgery, to
simplify the topology of a space:
\begin{figure}[ht]
\begin{center}
\begin{tabular}{ccc}
 \includegraphics[height=4.5cm]{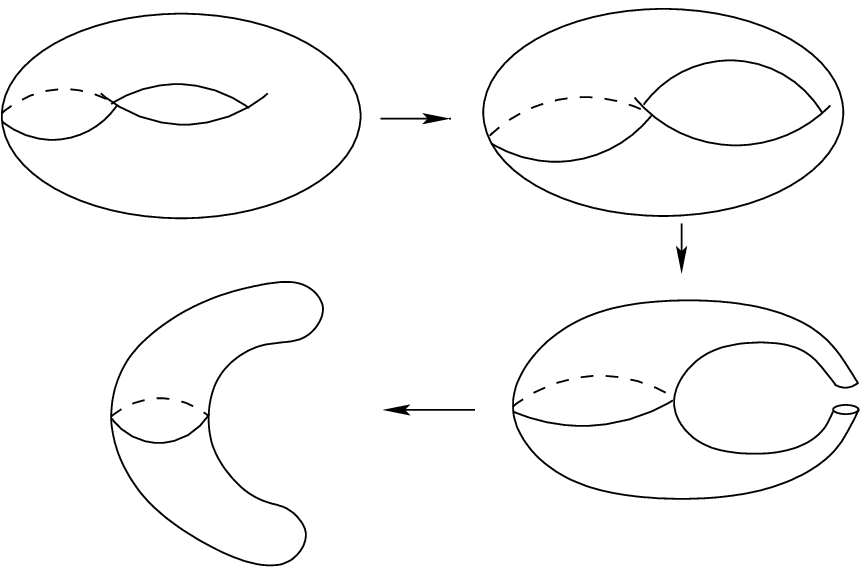}
\end{tabular}
\end{center}
\end{figure}

Two major ingredients were invented. One is called Dehn's lemma
which provides a tool to simplify any surface which cross itself
to one which does not.
\begin{figure}[ht]
\begin{center}
\begin{tabular}{ccc}
 \includegraphics[height=3.5cm]{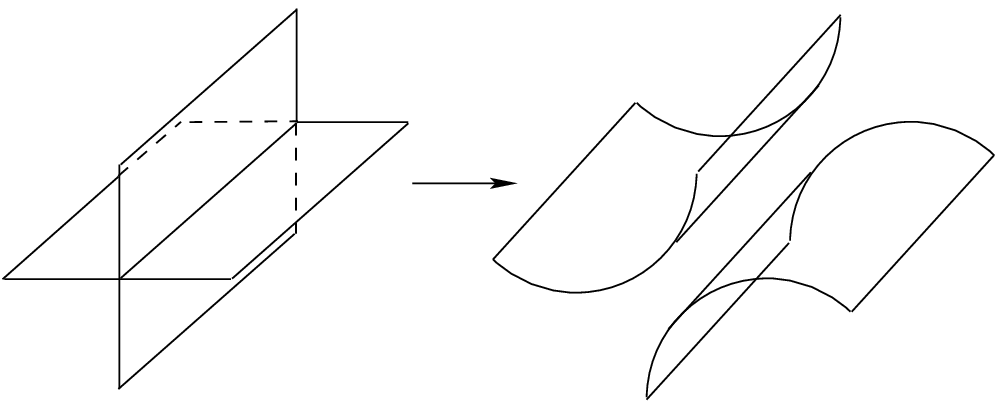}
\end{tabular}
\end{center}
\end{figure}

\noindent\textbf{Theorem (Dehn's lemma)} {\em If there exists a
map of a disk into a three dimensional space, which does not cross
itself on the boundary of the disk, then there exists another map
of the disc into the space which does not cross itself and is
identical to the original map on the boundary of the disc. }

This is a very subtle theorem, as it took almost fifty years until
Papakyriakopoulos came up with a correct proof after its
discovery.

The second tool is the construction of incompressible surfaces
introduced by Haken. It was used to cut three manifolds into
pieces. Walhausen proved important theorems by this procedure.
(Incompressible surfaces are embedded surfaces which have the
property whereby if a loop cannot be shrunk to a point on the
surface, then it cannot be shrunk to a point in the three
dimensional space, either.)

\section{Special Surfaces}

There are several important one dimensional and two dimensional
spaces that play important roles in understanding three dimensional
spaces.

\medskip

1. Circle

Seifert constructed many three dimensional spaces that can be
described as continuous family of circles. The above mentioned
phase space is an example of a \textbf{Seifert space}.

\medskip

2. Two dimensional spheres

We can build three dimensional spaces by removing balls from two
distinguished ones and gluing them along the boundary spheres.
\begin{figure}[ht]
\psfrag{S}{$S^2$}
\begin{center}
\begin{tabular}{ccc}
 \includegraphics[height=1.5cm]{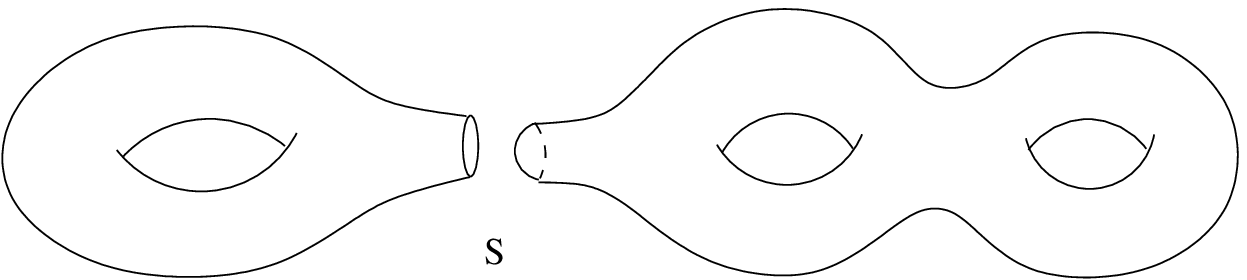}
\end{tabular}
\end{center}
\end{figure}
Conversely Kneser and Milnor proved that each three dimensional
space can be uniquely decomposed into irreducible components along
spheres. (A space is called irreducible if each embedded sphere is
the boundary of a three dimensional ball in this space.)

\medskip

3. Torus

A theorem of Jaco-Shalen, Johannson says that one can go one step
further by cutting a space along tori.
\begin{figure}[ht]
\psfrag{T}{$T^2$}
\begin{center}
\begin{tabular}{ccc}
 \includegraphics[height=2.0cm]{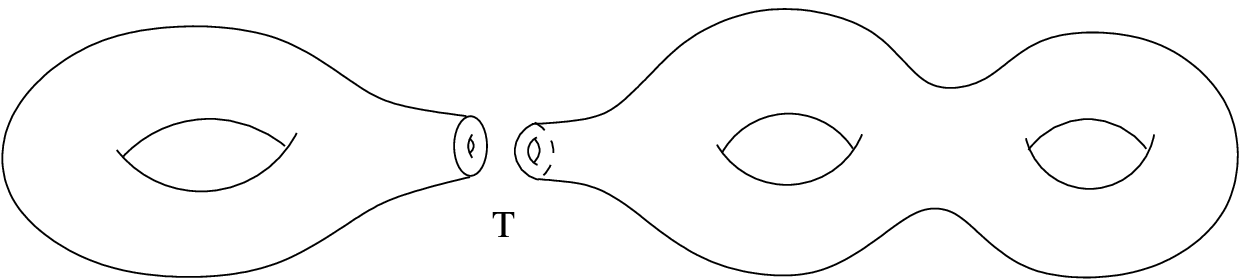}
\end{tabular}
\end{center}
\end{figure}

\section{Structure of Three Dimensional Spaces}

A very important breakthrough was made in the late 1970s by W.
Thurston. He make the following conjecture.

\textbf{Geometrization Conjecture (Thurston)}:  The structure of
three dimensional spaces is built on the following \textbf{atomic}
spaces:

(1) The Poincar\'e conjecture: three dimensional space where every
closed loop can be shrunk to a point; this space is conjectured to
be the three-sphere.

(2) The space-form problem: spaces obtained by identifying points
on the three-sphere. The identification is dictated by a finite
group of linear isometries which is similar to the symmetries of
crystals.

(3) Seifert spaces mentioned above and their quotients similar to
    (2).

(4) Hyperbolic spaces according to the conjecture of Thurston:
three-space whose boundaries may consist of tori such that every
two-sphere in the space is the boundary of a ball in the space and
each incompressible torus can be deformed to a boundary component;
it was conjectured to support a canonical metric with constant
negative curvature and it is obtained by identifying points on the
hyperbolic ball. The identification is dictated by a group of
symmetries of the ball similar to the symmetries of crystals.

An example of a space obtained by identifying points on the three
dimensional hyperbolic space

\begin{figure}[ht]
\begin{center}
\begin{tabular}{c}
\includegraphics[height=5.5cm]{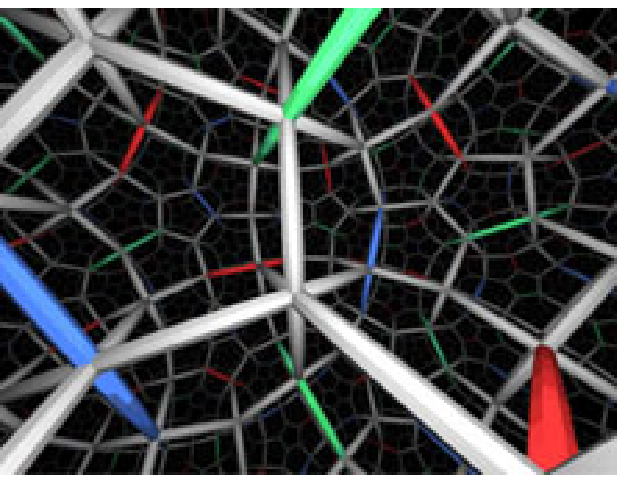}
\end{tabular}\\
Hyperbolic Space Tiled with Dodecahedra, \\by Charlie Gunn
(Geometry Center). \\ from the book "Three-dimensional geometry and topology" \\
by Thurston, Princeton University press
\end{center}
\end{figure}

Thurston's conjecture effectively reduced the classification of
three dimensional spaces to group theory, where many tools were
available. He and his followers proved the conjecture when the
three space is sufficiently large in the sense of Haken and
Walhausen. (A space is said to be sufficiently large if there is a
nontrivial incompressible surface embedded inside the space. Haken
and Walhausen proved substantial theorem for this class of
manifolds.) This theorem of Thurston covers a large class of three
dimensional hyperbolic manifolds.

However, as nontrivial incompressible surface is difficult to find
on a general space, the argument of Thurston is difficult to use
to prove the Poincar\'{e} conjecture.

\section{Geometric Analysis}

On the other hand, starting in the seventies, a group of geometers
applied nonlinear partial differential equations to build
geometric structures over a space. Yamabe considered the equation
to conformally deform metrics to metrics with constant scalar
curvature. However, in three dimension, metrics with negative
scalar curvature cannot detect the topology of spaces.

A noted advance was the construction of K\"{a}hler-Einstein
metrics on K\"{a}hler manifolds in 1976. In fact, I used such a
metric to prove the complex version of the Poincar\'e conjecture.
It is called the Severi conjecture in complex geometry. It says
that every complex surface that can be deformed to the complex
projective plane is itself the complex projective plane.

The subject of combining ideas from geometry and analysis to
understand geometry and topology is called \textbf{geometric
analysis}. While the subject can be traced back to 1950s, it has
been studied much more extensively in the last thirty years.

Geometric analysis is built on two pillars: nonlinear analysis and
geometry. Both of them became mature in the seventies based on the
efforts of many mathematicians. (See my survey paper
``Perspectives on geometric analysis" in {\it Survey in
Differential Geometry}, Vol. X. 2006.)

\section{Einstein metrics}

I shall now describe how ideas of geometric analysis are used to
solve the Poincar\'e conjecture and the geometrization conjecture
of Thurston.

In the case of a three dimensional space, we need to construct an
Einstein metric, a metric inspired by the Einstein equation of
gravity. Starting from an arbitrary metric on three space, we
would like to find a method to deform it to the one that satisfies
Einstein equation. Such a deformation has to depend on the
curvature of the metric.

Einstein's theory of relativity tells us that under the influence
of gravity, space-time must have curvature. Space moves
dynamically. The global topology of space changes according to the
distribution of curvature (gravity). Conversely, understanding of
global topology is extremely important and it provides constraints
on the distribution of gravity. In fact, the topology of space may
be considered as a source term for gravity.

From now on, we shall assume that our three dimensional space is
compact and has no boundary (i.e., closed).

In a three dimensional space, curvature of a space can be
different when measured from different directions. Such a
measurement is dictated by a quantity $R_{ij}$, called the
\textbf{Ricci tensor}. In general relativity, this gives rise to
the matter tensor of space.

An important quantity that is independent of directions is the
\textbf{scalar curvature} $R$. It is the trace of $R_{ij}$ and can
be considered as a way to measure the expansion or shrinking of
the volume of geodesic balls:
\[
   \textrm{Volume}(B(p,r)) \thicksim \frac{4\pi}{3} (r^3 - \frac{1}{30}R(p)
   r^5),
\]
where $B(p,r)$ is the ball of radius $r$ centered at a point $p$,
and $R(p)$ is the scalar curvature at $p$.

This can be illustrated by a dumbbell surface where, near the
neck, curvature is negative and where, on the two ends which are
convex, curvature is positive.
\begin{figure}[ht]
\begin{center}
\begin{tabular}{ccc}
 \includegraphics[width=3in]{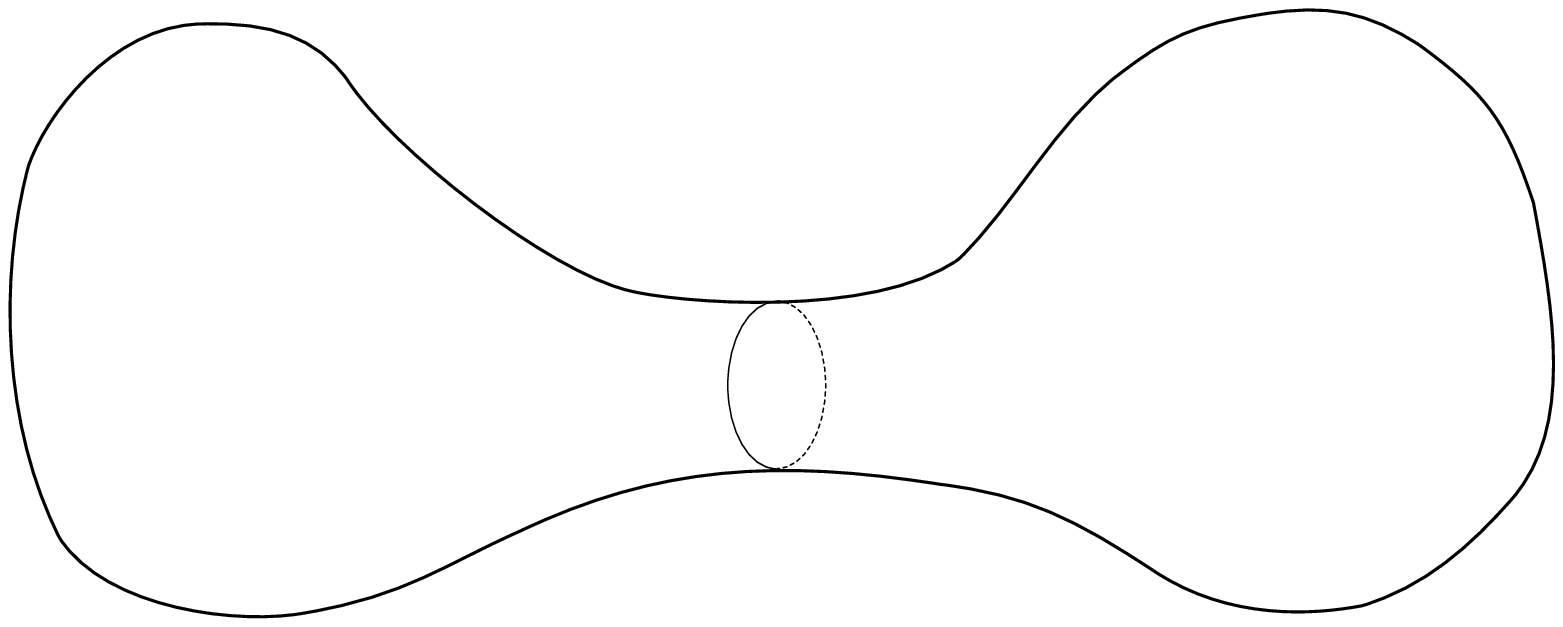}
\end{tabular}\\
Two-dimensional dumbbell surface
\end{center}
\end{figure}

Two-dimensional surfaces with negative curvature look like
saddles.
\begin{figure}[ht]
\begin{center}
\begin{tabular}{ccc}
 \includegraphics[height=3.0cm]{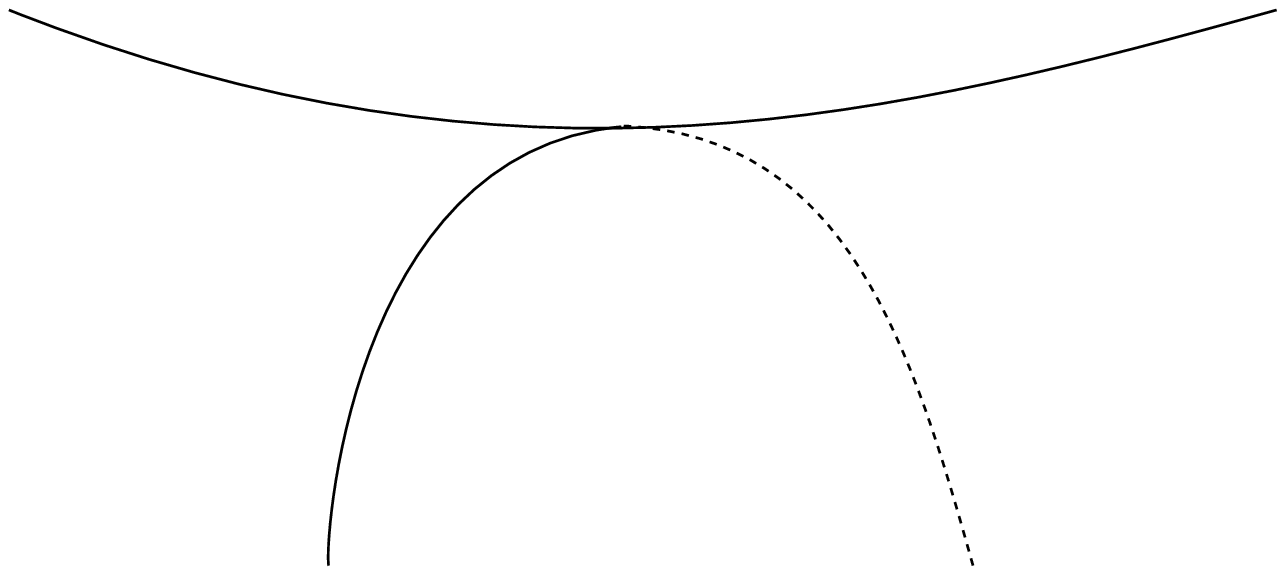}
\end{tabular}
\end{center}
\end{figure}
Hence a two dimensional neck has negatives curvature. However, in
three dimension, the slice of a neck can be a two dimensional
sphere with very large positive curvature. Since scalar curvature
is the sum of curvatures in all direction, the scalar curvature at
the three dimensional neck can be positive. This is an important
difference between a two-dimensional neck and a three-dimensional
neck.
\begin{figure}[ht]
\begin{center}
\begin{tabular}{ccc}
 \includegraphics[height=5.0cm]{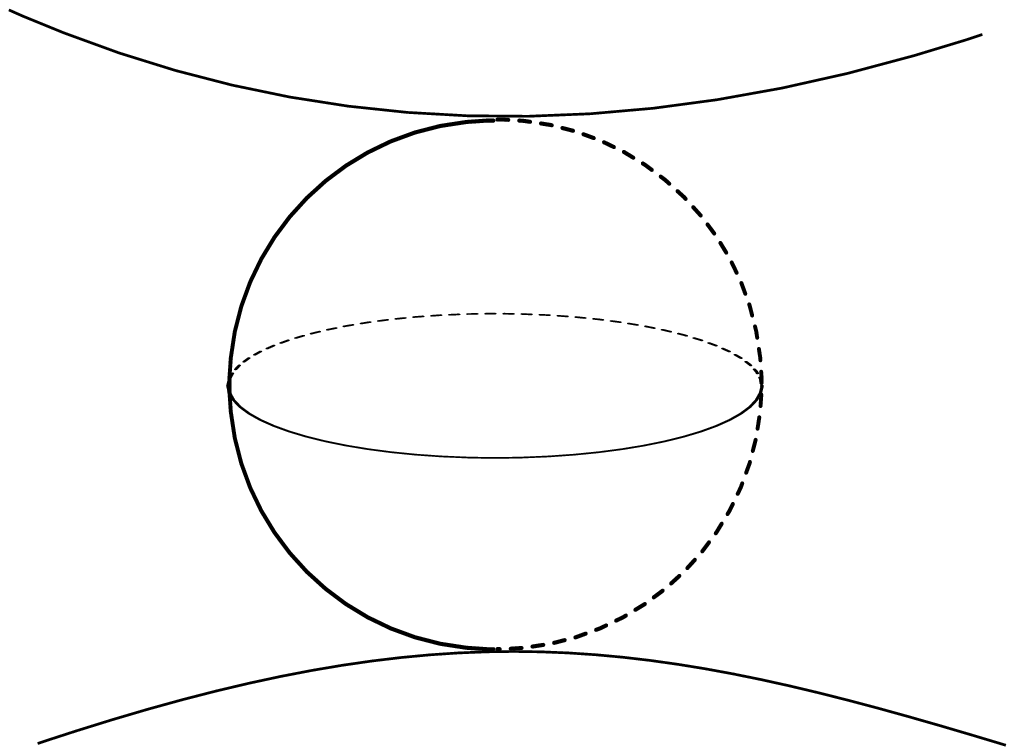}
\end{tabular}
\\
Three dimensional neck.
\end{center}
\end{figure}

\section{The dynamics of Einstein equation}

In general relativity, matter density consists of scalar curvature
plus the momentum density of space. The Dynamics of Einstein
equation drives space to form black holes which splits space into
two parts: the part where scalar curvature is positive and the
other part, where the space may have a black hole singularity and
is enclosed by the apparent horizon of the black hole, the
topology tends to support metrics with negative curvature.

There are two quantities in gravity that dictate the dynamics of
space: \textbf{metric} and \textbf{momentum}. Momentum is
difficult to control. Hence at this time, it is rather difficult
to use the Einstein equation of general relativity to study the
topology of spaces.

\section{Hamilton's Equation}

In 1979, Hamilton developed a new equation to study the dynamics
of space metric. The Hamilton equation is given by
\[
    \frac{\partial g_{ij}}{\partial t} = -2 R_{ij}.
\]

Instead of driving space metric by gravity, he drives it by its
Ricci curvature which is analogous to the heat diffusion.
Hamilton's equation therefore can be considered as a nonlinear
heat equation. Heat flows have a regularizing effect because they
disperse irregularity in a smooth manner.

Hamilton's equation was also considered by physicists. (It first
appeared in Friedan's thesis.) However, this point of view was
completely different. Physicists considered it as beta function
for deformations of the sigma model to conformal field theory.

\section{Singularity}

Despite the fact that Hamilton's equation tend to smooth out
metric structure, global topology and nonlinear terms in the
equation coming from curvature drive the space metric to points
where the space topology collapses. We call such points
\textbf{singularity of space}.

In 1982, Hamilton published his first paper on the equation.
Starting with a space with positive Ricci curvature, he proved
that under his equation, space, after dilating to keep constant
volume, never encounters any singularity and settles down to a
space where curvature is constant in every direction.

Such a space must be either a $3$-sphere or a space obtained by
identifying the sphere by some finite group of isometries.

After seeing the theorem of Hamilton, I was convinced that
Hamilton's equation is the right equation to carry out the
geometrization program. (This paper of Hamilton is immediately
followed by the paper of Huisken on deformation of convex surfaces
by mean curvature. The equation of mean curvature flow has been a
good model for understanding Hamilton's equation.)

We propose to deform any metric on a three dimensional space which
shall break up the space eventually. It should lead to the
topological decomposition according to Kneser, Milnor,
Jacob-Shalen and Johannson. The asymptotic state of Hamilton's
equation is expected to be broken up into pieces which will either
collapse or produce metrics which satisfy the Einstein equation.

In three dimensional spaces, Einstein metrics are metrics with
constant curvature. However, along the way, the deformation will
encounter singularities. The major question is how to find a way
to describe all possible singularities. We shall describe these
spectacular developments.

\section{Hamilton's Program}

Hamilton's idea is to perform surgery to cut off the singularities
and continue his flow after the surgery. If the flow develops
singularities again, one repeats the process of performing surgery
and continuing the flow.

If one can prove there are only a finite number of surgeries in
any finite time interval, and if the long-time behavior of
solutions of the Hamilton's flow with surgery is well understood,
then one would be able to recognize the topological structure of
the initial manifold. Thus Hamilton's program, when carried out
successfully, will lead to a proof of the Poincar\'e conjecture
and Thurston's geometrization conjecture.

The importance and originality of Hamilton's contribution can
hardly be exaggerated. In fact, Perelman said:

\textbf{``The implementation of Hamilton's program would imply the
geometrization conjecture for closed three-manifolds."}

\textbf{``In this paper we carry out some details of Hamilton's
program".}

We shall now describe the chronological development of Hamilton's
program. There were several stages:

\medskip

\begin{center}
\textbf{I. A Priori Estimates}
\end{center}

In the early 1990s, Hamilton systematically developed methods to
understand the structure of singularities. Based on my suggestion,
he proved the fundamental estimate (the Li-Yau-Hamilton estimate)
for his flow when curvature is nonnegative. The estimate provides
a priori control of the behavior of the flow.

An a prior estimate is the key to proving any existence theorem
for nonlinear partial differential equations. An intuitive example
can be explained as follows: when a missile engineer designs
trajectory of a missile, he needs to know what is the most likely
position and velocity of the missile after ten seconds of its
launch. Yet a change in the wind will cause reality to differ from
his estimate. But as long as the estimate is within a range of
accuracy, he will know how to design the missile. How to estimate
this range of accuracy is called a prior estimate.

\medskip

\textbf{The Li-Yau-Hamilton Estimate}

In proving existence of a nonlinear differential equation, we need
to find an a priori estimate for some quantity which governs the
equation. In the case of Hamilton's equation, the important
quantity is the scalar curvature $R$. An absolute bound on the
curvature gives control over the nonsingularity of the space. On
the other hand, the relative strength of the scalar curvature
holds the key to understand the singularity of the flow. This is
provided by the Li-Yau-Hamilton estimate:

For any one-form $V_a$ we have
$$
\frac{\partial R}{\partial t}+\frac{R}{t} +2\nabla_aR\cdot
V_a+2R_{ab}V_aV_b\geq 0.
$$ In particular, $tR(x,t)$ is pointwise nondecreasing in time.

In the process of applying such an estimate to study the structure
of singularities, Hamilton discovered (also independently by Ivey)
a curvature pinching estimate for his equation on
three-dimensional spaces. It allows him to conclude that a
neighborhood of the singularity looks like space with nonnegative
curvature. For such a neighborhood, the Li-Yau-Hamilton estimate
can be applied.

Then, under an additional non-collapsing condition, Hamilton
described the structure of all possible singularities. However, he
was not able to show that all these possibilities actually occur.
Of particular concern to him was a singularity which he called the
cigar.

\medskip

\begin{center}
\textbf{II. Hamilton's works on Geometrization}
\end{center}

In 1995, Hamilton developed the procedure of geometric surgery
using a foliation by surfaces of constant mean curvature, to study
the topology of four-manifolds of positive isotropic curvature.

In 1996, he went ahead to analyze the global structure of the
space time structure of his flow under suitable regularity
assumptions (he called them nonsingular solutions). In particular,
he showed how three-dimensional spaces admitting a nonsingular
solution of his equation can be broken into pieces according to
the geometrization conjecture.

These spectacular works are based on deep analysis of geometry and
nonlinear differential equations. Hamilton's two papers provided
convincing evidence that the geometrization program could be
carried out using his approach.

\medskip

\textbf{Main Ingredients of these works of Hamilton}

In this deep analysis he needed several important ingredients:

(1)  a compactness theorem on the convergence of metrics developed
by him, based on the injectivity radius estimate proved by
Cheng-Li-Yau in 1981. (The injectivity radius at a point is the
radius of the largest ball centered at that point that the ball
would not collapse topologically.)

(2) a beautiful quantitative generalization of Mostow's rigidity
theorem which says that there is at most one metric with constant
negative curvature on a three-dimensional space with finite
volume. This rigidity theorem of Mostow is not true for two
dimensional surfaces.

(3) In the process of breaking up the space along the tori, he
needs to prove that the tori are incompressible. The ingredients
of his proof depend on the theory of minimal surfaces as was
developed by Meeks-Yau and Schoen-Yau.

At this stage, it seems clear to me that Hamilton's program for
the Poincar\'e and geometrization conjectures could be carried
out. The major remaining obstacle was to obtain certain
injectivity radius control, in terms of local curvature bound, in
order to understand the structure of the singularity and the
process of surgery to remove the singularity. Hamilton and I
worked together on removing this obstacle for some time.

\medskip

\begin{center}
\textbf{III. Perelman's Breakthrough}
\end{center}

In November of 2002, Perelman put out a preprint, ``The entropy
formula for Hamilton's equation and its geometric applications",
wherein major ideas were introduced to implement Hamilton's
program.

Parallel to what Li-Yau did in 1986, Perelman introduced a
space-time distance function obtained by path integral and used it
to verify the noncollapsing condition in general. In particular,
he demonstrated that cigar type singularity does not exist in
Hamilton's equation.

His distance function can be described as follows.

Let $\sigma$ be any space-time path joining $p$ to $q$, we define
the action to be
$$ \int^{\tau}_0 \sqrt{s}(R+| \dot{\sigma}(s)|^2)ds.$$ By
minimizing among all such paths joining $p$ to $q$, we obtain
$L(q, \tau)$.

Then Perelman defined his reduced volume to be $$\int
(4\pi\tau)^{-\frac{n}{2}}\exp
\left\{-\frac{1}{2\sqrt{\tau}}L(q,\tau)\right\}$$ and observed
that under the Hamilton's equation it is nonincreasing in $\tau$.
In this proof Perelman used the idea in the second part of
Li-Yau's paper in 1986. As recognized by Perelman: ``in Li-Yau,
where they use `length', associated to a linear parabolic
equation, is pretty much the same as in our case."

\medskip

\textbf{Rescaling Argument}

Furthermore, Perelman developed an important refined rescaling
argument to complete the classification of Hamilton on the structure
of singularities of Hamilton's equation and obtained a uniform and global
version of the structure theorem of singularities.

\begin{figure}[ht]
\begin{center}
\begin{picture}(320, 240)
\qbezier(15,90)(0,50)(20,20)\qbezier(20,20)(50,10)(100,50)\qbezier(100,50)(130,60)(170,70)
\qbezier(15,90)(50,140)(100,130)\qbezier(100,130)(130,140)(170,150)\qbezier(130,130)(110,90)(130,70)
\qbezier(130,140)(126,135)(130,130)\qbezier(130,70)(126,65)(130,60)
\qbezier[8](130,140)(134,135)(130,130)\qbezier[8](130,70)(134,65)(130,60)
\qbezier(130,130)(150,135)(170,150)\qbezier(130,70)(150,70)(170,70)
\qbezier(170,70)(190,78)(200,80)\qbezier(170,70)(190,68)(200,70)
\qbezier(200,80)(210,78)(230,70)\qbezier(200,70)(210,68)(230,70)
\qbezier(200,80)(196,75)(200,70)\qbezier[8](200,80)(204,75)(200,70)
\qbezier(230,70)(260,70)(280,70)\qbezier(230,70)(260,60)(280,50)\qbezier(280,70)(310,60)(280,50)
\qbezier(260,70)(256,65)(260,60)\qbezier[8](260,70)(264,65)(260,60)
\qbezier(170,150)(200,160)(220,170)\qbezier(170,150)(200,150)(220,150)
\qbezier[8](200,160)(204,155)(200,150)\qbezier(200,160)(196,155)(200,150)
\qbezier(220,170)(260,230)(290,200)\qbezier(220,150)(260,150)(290,170)\qbezier(290,200)(310,200)(290,170)
\qbezier(230,160)(260,160)(290,190)\qbezier(250,165)(260,185)(280,180)
\qbezier(30,60)(50,60)(80,100)\qbezier(40,65)(50,85)(75,90)
\qbezier(238,160)(230,155)(235,150)\qbezier[8](238,160)(241,155)(235,150)
\qbezier(280,180)(290,175)(290,170)\qbezier[10](280,180)(278,175)(290,170)

\put(30,140){\makebox(0,0)[bl]{$\Omega_{\rho}$}}
\put(40,140){\vector(1,-1){10}}
\put(200,190){\makebox(0,0)[bl]{$\Omega_{\rho}$}}
\put(210,190){\vector(1,-1){10}}
\put(160,115){\makebox(0,0)[bl]{$\varepsilon$-horn}}
\put(160,125){\vector(-1,1){10}}
\put(270,130){\makebox(0,0)[bl]{$\varepsilon$-tube}}
\put(270,140){\vector(-1,1){10}}
\put(160,48){\makebox(0,0)[bl]{double $\varepsilon$-horn}}
\put(190,58){\vector(0,1){10}}
 \put(250,30){\makebox(0,0)[bl]{capped $\varepsilon$-horn}}
 \put(270,40){\vector(0,1){10}}
\end{picture}\\
The Structure of Singularity
\end{center}
\end{figure}

\medskip

\textbf{Hamilton's Geometric Surgery}

 Now we need to find a way to perform geometric surgery. In
1995, Hamilton had already initiated a surgery procedure for his
equation on four-dimensional spaces and presented a concrete
method for performing such surgery.

One can see that Hamilton's geometric surgery method also works
for Hamilton's equation on three-dimensional spaces. However, in
order for surgeries to be done successfully, a more refined
technique is needed.
\begin{figure}[ht]
\begin{center}
\begin{picture}(320,120)
\qbezier(40,10)(60,30)(310,50) \qbezier(40,90)(60,70)(310,50)
\qbezier(185,39)(190,39)(190,50)\qbezier(185,61)(190,61)(190,50)
\qbezier[5](185,39)(181,39)(181,50)\qbezier[5](185,61)(181,61)(181,50)
\qbezier(80,24)(93,27)(93,50)\qbezier(80,76)(93,73)(93,50)
\qbezier[15](80,24)(67,27)(67,50)\qbezier[15](80,76)(67,73)(67,50)
\qbezier(135,33)(143,33)(143,50)\qbezier(135,67)(143,67)(143,50)
\qbezier[8](135,33)(127,33)(127,50)\qbezier[8](135,67)(127,67)(127,50)
\put(240,20){\makebox(0,0)[bl]{$\varepsilon$-horn}}
\put(250,30){\vector(-1,1){10}}
\put(130,10){\makebox(0,0)[bl]{neck}}
\put(140,20){\vector(-1,1){10}}
\end{picture}
\begin{picture}(320,120)
\qbezier(40,10)(60,23)(135,33)\qbezier[50](135,33)(185,39)(310,50)
\qbezier(40,90)(60,77)(135,67)\qbezier[50](135,67)(185,61)(310,50)
\qbezier[5](185,39)(190,39)(190,50)\qbezier[5](185,61)(190,61)(190,50)
\qbezier[5](185,39)(181,39)(181,50)\qbezier[5](185,61)(181,61)(181,50)
\qbezier(80,24)(93,27)(93,50)\qbezier(80,76)(93,73)(93,50)
\qbezier[15](80,24)(67,27)(67,50)\qbezier[15](80,76)(67,73)(67,50)
\qbezier(135,33)(143,33)(143,50)\qbezier(135,67)(143,67)(143,50)
\qbezier[8](135,33)(127,33)(127,50)\qbezier[8](135,67)(127,67)(127,50)
\qbezier(135,33)(175,38)(175,50)\qbezier(135,67)(175,62)(175,50)
\put(155,20){\makebox(0,0)[bl]{the gluing cap}}
\put(170,30){\vector(-1,1){10}}
\end{picture}\\
geometric surgery
\end{center}
\end{figure}

\medskip

\textbf{Discreteness of Surgery Times}

The challenge is to prove that there are only a finite number of
surgeries on each finite time interval. The problem is that, when
one performs the surgeries with a given accuracy at each surgery
time, it is possible that the error may add up to so fast that
they force the surgery times to accumulate.

\medskip

\textbf{Rescaling Arguments}

In March of 2003, Perelman put out another preprint, titled
``Ricci flow with surgery on three manifolds", where he designed
an improved version of Hamilton's geometric surgery procedure so
that, as time goes on, successive surgeries are performed with
increasing accuracy.

Perelman introduced a rescaling argument to prevent the surgery
time from accumulating.

When using the rescaling argument for surgically modified
solutions of Hamilton's equation, one encounters the difficulty of
applying Hamilton's compactness theorem, which works only for
smooth solutions.

The idea of overcoming this difficulty consists of two parts:

1. (Perelman): choose the cutoff radius in the neck-like regions
small enough to push the surgical regions far away in space.

2. (Cao-Zhu): establish results for the surgically modified
solutions so that Hamilton's compactness theorem is still
applicable. To do so, they need a deep understanding of the
prolongation of the surgical regions, which in turn relies on the
uniqueness theorem of Chen-Zhu for solutions of Hamilton's
equation on noncompact manifolds.

\medskip

\textbf{Conclusion of the proof of the Poincar\'e Conjecture}

One can now prove Poincar\'{e} conjecture for simply connected
three dimensional space, by combining the discreteness of
surgeries with finite time extinction result of Colding-Minicozzi
(2005).

\medskip

\begin{center}
\textbf{IV. Proof of the geometrization conjecture: Thick-thin
Decomposition}
\end{center}

To approach the structure theorem for general spaces, one still
needs to analyze the long-time behavior of surgically modified
solutions to Hamilton's equation. As mentioned in II, Hamilton
studied the long time behavior of his equation for a special class
of (smooth) solutions -- nonsingular solutions.

\begin{figure}[ht]
\begin{center}
\begin{picture}
 (320,270)
\qbezier(50,210)(20,210)(5,170)\qbezier(5,170)(0,95)(40,90)\qbezier(40,90)(180,40)(230,15)
\qbezier(230,15)(270,10)(280,25)\qbezier(280,25)(280,50)(250,50)\qbezier(250,50)(180,20)(110,90)
\qbezier(110,90)(70,125)(110,160)\qbezier(110,160)(190,220)(280,190)\qbezier(280,190)(310,190)(315,210)
\qbezier(315,210)(315,260)(270,240)\qbezier(270,240)(190,170)(50,210)
\qbezier(30,140)(50,140)(70,165)\qbezier(40,145)(50,150)(55,152)
\qbezier(230,30)(250,25)(270,30)\qbezier(240,30)(250,35)(260,30)
\qbezier(260,220)(280,210)(290,220)\qbezier(270,215)(280,220)(285,220)
\qbezier(120,60)(127,68)(125,75)\qbezier[8](120,60)(117,73)(125,75)
\qbezier(195,30)(200,35)(200,40)\qbezier[5](195,30)(192,36)(200,40)
\qbezier(147,195)(150,190)(147,183)\qbezier[5](147,195)(144,190)(147,183)
\qbezier(226,210)(230,205)(230,200)\qbezier[5](226,210)(222,205)(230,200)
\put(270,180){\vector(0,1){10}}\put(250,170){\makebox(0,0)[bl]{thick
part}}\put(235,160){\makebox(0,0)[bl]{(hyperbolic piece)}}
\put(250,60){\vector(0,-1){10}}\put(230,70){\makebox(0,0)[bl]{thick
part}}\put(215,60){\makebox(0,0)[bl]{(hyperbolic piece)}}
\put(40,80){\vector(1,1){10}}\put(20,70){\makebox(0,0)[bl]{thick
part}}\put(5,60){\makebox(0,0)[bl]{(hyperbolic piece)}}
\put(140,40){\vector(1,1){10}}\put(120,30){\makebox(0,0)[bl]{thin
part}}
\put(170,210){\vector(1,-1){10}}\put(150,210){\makebox(0,0)[bl]{thin
part}}
\end{picture}\\
Thick-thin decomposition
\end{center}
\end{figure}

In 1996, Hamilton proved that any three-dimensional nonsingular
solution admits of a thick-thin decomposition where the thick part
consists of a finite number of hyperbolic pieces and the thin part
collapses. Moreover, by adapting Schoen-Yau's minimal surface
arguments, Hamilton showed that the boundary of hyperbolic pieces
are incompressible tori. Consequently, any nonsingular solution is
geometrizable.

Even though the nonsingularity assumption seems restrictive, the
ideas and arguments of Hamilton are used in an essential way by
Perelman to analyze the long-time behavior for general surgical
solutions. In particular, he also studied the thick-thin
decomposition.

For the thick part, based on the Li-Yau-Hamilton estimate,
Perelman established a crucial elliptic type estimate, which
allowed him to conclude that the thick part consists of hyperbolic
pieces. For the thin part, since he could only obtain a lower
bound on the sectional curvature, he needs a new collapsing
result. Assuming this new collapsing result, Perelman claimed that
the solutions to Hamilton's equation with surgery have the same
long-time behavior as nonsingular solutions in Hamilton's work, a
conclusion which would imply the validity of Thurston's
geometrization conjecture.

Although the proof of this new collapsing result was promised by
Perelman, it still has yet to appear. (Shioya-Yamaguchi has
published a proof of the collapsing result in the special case
when the space is closed.) Nonetheless, based on the previous
results, Cao-Zhu gave a complete proof of Thurston's
geometrization conjecture.

\section*{
    \textbf{Conclusion}
}

The success of Hamilton's program is the culmination of efforts by
geometric analysts in the past thirty years. It should be
considered as the crowning achievement of the subject of geometric
analysis, a subject that is capable of proving hard and difficult
topological theorems by geometry and analysis solely.

Hamilton's equation is a complicated nonlinear system of partial
differential equations. This is the first time that mathematicians
have been able to understand the structure of singularity and
development of such a complicated system.

Similar systems appear throughout the natural world. The methods
developed in the study of Hamilton equation should shed light on
many natural systems such as the Navier-Stokes equation and the
Einstein equation.

In addition, the numerical implementation of the Hamilton flow
should be useful in computer graphics, as was demonstrated by
Gu-Wang-Yau for two dimensional figures.

\section*{
\textbf{Impact on the future of geometry} }

\textbf{Poincar\'e:}

``Thought is only a flash in the middle of a long night, but the
flash that means everything.''

The Flash of Poincar\'e in 1904 has illuminated a major portion of
the topological developments in the last century.

Poincar\'e also initiated development of the theory of Riemann
surfaces. It has been one of the major pillars of all mathematics
development in the twentieth century. I believe that the full
understanding of the three dimensional manifolds will play a
similar role in the twenty-first century.

\medskip

\begin{center}
\textbf{Remark}
\end{center}

In Perelman's work, many key ideas of the proofs are sketched or
outlined, but complete details of the proofs are often missing.
The recent paper of Cao-Zhu, submitted to The Asian Journal of
Mathematics in 2005, gives the first complete and detailed account
of the proof of the Poincar\'e conjecture and the geometrization
conjecture. They substituted several arguments of Perelman with
new approaches based on their own studies. The materials were
presented by Zhu in a Harvard seminar from September 2005 to March
2006, where faculties and postdoctoral fellows of Harvard
University and MIT attended regularly. Some of the key arguments,
that has been important for the completion of the Poincare
conjecture, has already appeared in the paper of Chen-Zhu
\cite{Chen-Zhu}.

In the last three years, many mathematicians have attempted to see
whether the ideas of Hamilton and Perelman can hold together.
Kleiner and Lott (in 2004) posted on their web page some notes on
several parts of Perelman's work. However, these notes were far
from complete. After the work of Cao-Zhu was accepted and
announced by the journal in April, 2006 (it was distributed on
June 1, 2006). On May 24, 2006, Kleiner and Lott put up another,
more complete, version of their notes. Their approach is different
from Cao-Zhu's. It will take some time to understand their notes
which seem to be sketchy at several important points. Most
recently, a manuscript of Morgan-Tian appeared in the web. In a
letter to the author, Jim Carlson of the Clay institute stated
that the first version of this manuscript was submitted to the
Clay institute on May 19, 2006, and the revised version was
submitted on July 23, 2006.

\end{document}